\documentclass[11pt]{article}

\usepackage{amsfonts,amsmath,graphicx}
\graphicspath{{.}}
   \makeatletter
   \let\accent@spacefactor\relax
   \makeatother
\usepackage{amsfonts}
\usepackage{amsmath}   
\usepackage{amsthm}    
\usepackage{amscd}     
\usepackage{euscript}

\def\C{{\Bbb C}}
\def\P{{\Bbb P}}

\newtheorem{defi}{D\'{e}finition}[section]
\newtheorem{pro}[defi]{Proposition}

\newtheorem{theo}[defi]{Theorem}

\newtheorem{coro}[defi]{Corollary}

\title{
Invariant multidimensional matrices}
\author{Jean Vall\`es}
\date{}

\begin{document}

\maketitle
{\small {\it ``On se persuade mieux, pour l'ordinaire , par les
raisons qu'on a soi-m\^eme trouv\'ees que par celles qui sont venues dans l'esprit des
autres''(Pascal).}}

\bigskip

{\small {\bf Abstract :} In \cite{AO} the authors study Steiner bundles via their
unstable hyperplanes and proved that (see \cite{AO}, Tmm 5.9) :\\
{\it A rank $n$ Steiner bundle on $\P^n$ which is $SL(2,\C)$ invariant is a
Schwarzenberger bundle.} \\
In this note we give  a very short proof  of this result based on Clebsch-Gordon
problem for $SL(2,\C)$-modules.}
\section{Introduction}
Let $A$, $B$, and $C$ three vector spaces over $\C$ and $\phi : A\otimes B
\rightarrow C^{*}$ a linear surjective map. We consider the  sheaf ${\cal S}_{\phi}$
on 
$\P(A)\times \P(B)$ defined by,  
$$ \begin{CD}
0 @>>> {\cal S}_{\phi} @>>>C\otimes O_{\P(A)\times \P(B)}
@>^t\phi>> O_{\P(A)\times \P(B)}(1,1)
\end{CD}
$$
with fibers ${\cal S}_{\phi}(a\otimes b)=\{c\in C\mid
\phi(a\otimes b)(c)=0\}$.

\smallskip
 
{\bf Remark 1.} If 
${\rm dim}_{\C}C<{\rm dim}_{\C} A+{\rm dim}_{\C} B-1$ then  
${\rm ker}(\phi)$ meets the set of decomposable tensors, so there
exist 
$a\in A$ , $a\neq 0$ and $b\in B$, $b\neq 0$ such that 
$\phi(a\otimes b)=0$.

\smallskip

{\bf Remark 2.} When ${\cal S}_{\phi}$ is a vector bundle 
 it gives two ``associated'' Steiner bundles $S_A$ on $\P(A)$ and $S_B$
on $\P(B)$ after projections (see \cite{DK}, prop. 3.20). 

\medskip

We denote by $ (\P(A)\times \P(B)\times\P(C))^{\vee}$ the variety of hyperplanes
tangent to the Segre $\P(A)\times \P(B)\times\P(C)$ and 
by $(\{a\}\times
\{b\}\times\P(C))^{\vee}$ the set of hyperplanes in
$\P(A\otimes B\otimes C)$ containing $\{a\}\times \{b\}\times\P(C)$.

\smallskip

The next proposition is a reformulation  of many results from, \cite{GKZ} (see for
instance Thm 3.1' page 458, prop 1.1 page 445),  \cite{AO} (see thm page
1) and \cite{DO} (see cor.3.3).
\begin{pro} \label{hyper}
Let $A$, $B$, and $C$ three vector spaces over $\C$ with 
 ${\rm dim}_{\C} A=n+1$, ${\rm dim}_{\C} B=m+1$ and 
${\rm dim}_{\C} C\ge n+m+1$ and $\phi : A\otimes B\rightarrow C^{*}$ a surjective
linear map. Then the following propositions are equivalent :\\
1) ${\cal S}_{\phi}$ is a vector bundle over $\P(A)\times \P(B)$.\\   
2) $\phi(a\otimes b)\neq 0$ for all
$a\in A$, $a\neq 0$ and $b\in B$, $b\neq 0$.\\
3) $\phi \notin (\{a\}\times
\{b\}\times\P(C))^{\vee}$ for all
$a\in A$, $a\neq 0$ and $b\in B$, $b\neq 0$.\\
4)  $\phi \notin (\P(A)\times \P(B)\times\P(C))^{\vee}$
\end{pro}
{\bf Remark 3.} When  ${\rm dim}_{\C}C={\rm dim}_{\C} A+{\rm dim}_{\C} B-1$ the
variety   $(\P(A)\times
\P(B)\times\P(C))^{\vee}$ is an hypersurface in $\P((A\otimes B\otimes C)^{\vee})$.
This hypersurface  is defined by the vanishing of the hyperdeterminant, say
$Det(\Phi)$ where
$\Phi$ is the generic tridimensional matrix (see \cite{GKZ}, chapter 1 and 14). 

\smallskip
 
{\bf Proof.} It is clear that $1)$ $2)$ and $3)$ are equivalent. It remains to
show that $3)$ and $4)$ are equivalent too.
\\
Since
$\partial (abc)=(\partial (a))bc+a(\partial(b))c+ab(\partial (c))
$ an hyperplane $H$ is tangent to the Segre in a point $(a,b,c)$ if and only if it
contains 
$\P(A)\times \{b\}\times \{c\}$ and $\{a\}\times \P( B)\times \{c\}$ and $\{a\}\times
\{b\}\times\P(C)$. We prove here that the third condition implies the two others. Let $H$
an hyperplane containing $\{a\}\times \{b\}\times\P(C)$, we show that there exists $c\in
C$ such that $H$ contains 
$\P(A)\times \{b\}\times \{c\}$ and $\{a\}\times \P( B)\times \{c\}$. Let $\phi$ the
trilinear application corresponding to $H$. Since $\phi(a\otimes b)(C)=0$ we have a 
${\rm dim}_{\C}C$-dimensional family of bilinear forms vanishing on $(a,b)$. Now
finding a bilinear form of the above family  (i.e. finding $c\in C$) which verify 
 $\phi(A\otimes b)(c)=0$ and  $\phi(a\otimes B)(c)=0$ imposes at most $n+m$ conditions.
Since ${\rm dim}_{\C} C\ge n+m+1$, this point $c$ exists.~$\Box$
\section{Invariant tridimensional matrix under $SL(2,\C)$-action}
In the second part of this note we will consider the boundary case
$${\rm dim}_{\C}C={\rm dim}_{\C} A+{\rm dim}_{\C} B-1$$ Then, instead of
writing $\phi $ induces a vector bundle or $\phi \notin  
(\P(A)\times
\P(B)\times\P(C))^{\vee}$ we will write equivalently $Det(\phi)\neq 0$.

\smallskip

We denote by $S_i$ the irreducible $SL(2,\C)$-representations of degree $i$ and by 
$(x^{i-k}y^k)_{k=0,\cdots,i}$ a basis of $S_i$.
\begin{theo}
Let $A$, $B$ and $C$ be three non trivial $SL(2,\C)$-modules with dimension $n+1$, $m+1$
and $n+m+1$ and $\phi \in \P(A\otimes B \otimes C)$ an invariant hyperplane under
$SL(2,\C)$. Then, 
$$ Det(\phi)\neq 0 \Leftrightarrow \phi \,\,{\rm is }\,\,{\rm the}\,\, {\rm
multiplication}\,\, S_n\otimes S_m \rightarrow S_{n+m} $$
\end{theo}
{\bf Proof.} When $\phi \in \P(S_n\otimes S_m\otimes S_{n+m})$ is just
the multiplication $S_n\otimes S_m \rightarrow S_{n+m} $  it is well known that it
corresponds to Schwarzenberger bundles (see \cite{DK}, prop 6.3).

\smallskip

Conversely, let $A=\oplus_{i\in I} S_i\otimes U_i$, $B=\oplus_{j\in J} S_j\otimes
V_j$ where $U_i$, $V_j$ are trivial $SL(2,\C)$-representations of dimension
$n_i$ and $m_j$. Let   
$x^{i} \in S_i$,
$x^{j}\in S_j$ be two highest
weight vectors and $u\in U_i$, $v\in V_j$. Since $Det(\phi)\neq 0$, $\phi
((x^{i}\otimes u)\otimes (x^{j}\otimes v))\neq 0 $
  and by $SL(2,\C)$-invariance
$\phi ((x^{i}\otimes u)\otimes (x^{j}\otimes v))=x^{i+j}\phi (u\otimes v)\in
S_{i+j}\otimes W_{i+j}$. 
By hypothesis $\phi (u\otimes v)\neq 0$ for all $u\in U_i$ and $v\in V_j$ so, 
by the Remark 1,  it
implies that ${\rm dim }W_{i+j}\ge n_i+m_j-1$, and $S_{i+j}^{n_i+m_j-1}\subset
C^{*}$. 

\smallskip

Assume now that $B$ contains 		at least two
distinct irreducible representations. Let $i_{0}$ and $j_{0}$ the greatest integers in $I$
and $J$. We consider the submodule $B_1$ such that $B_1\oplus S_{j_0}^{m_{j_0}}=B$. Then
the restricted map 
$A\otimes B_1 \rightarrow C^{*}$ is not surjective because the image is concentrated in 
the submodule $C_1^{*} $ of $C^{*}$ defined by 
$C_1^{*}\oplus S_{i_0+j_0}^{n_{i_0}+m_{j_0}-1}=C^{*}$. Now since 
$${\rm dim}_{\C}(C_1)< {\rm dim}_{\C}(A)+{\rm dim}_{\C}(B_1)-1$$ there exist $a\in
A$,
$b\in B_1\subset B$ such that $\phi (a\otimes b)=0$. A contradiction with the
hypothesis
$Det(\phi)\neq 0$.

\smallskip

So $A=S_i^{n_i}$,  $B=S_j^{m_j}$ and $S_{i+j}^{n_i+m_j-1}\subset
C^{*}$. Since  ${\rm dim}_{\C}C={\rm dim}_{\C} A+{\rm dim}_{\C} B-1$, we have 
 $(i+1)n_i+(j+1)m_j -1={\rm dim}_{\C}C \ge (i+j+1)(n_i+m_j-1)$ which is possible if
and only if $n_i=m_j=1$ and $C=S_{i+j}$. ~$\Box$
\begin{coro}
 A rank $n$ Steiner bundle on $\P^n$ which is $SL(2,\C)$ invariant is a
Schwarzenberger bundle.
\end{coro}
{\it Proof.} Let $S$ a rank $n$ Steiner bundle on $\P^n$, i.e $S$ appears in an
exact sequence
$$ \begin{CD}
0 @>>> S @>>>C\otimes O_{\P(A)}
@>>> B^{*}\otimes O_{\P(A)}(1)
@>>>0
\end{CD}
$$
where $\P(A)=\P^n$, $\P(B)=\P^m$ and $\P (C)=\P^{n+m}$. 
If $SL(2,\C)$ acts on $S$ the vector spaces $A$, $B$ and $C$ are $SL(2,\C)$-modules
since $A$ is the basis, $B^{*}=H^1S(-1)$ and $C^{*}=H^0(S^{*})$. If $S$ is
$SL(2,\C)$-invariant the linear surjective map
$$ A\otimes (H^1S(-1))^{*} \rightarrow H^0(S^{*}) $$ is $SL(2,\C)$-invariant
too.~$\Box$

\medskip

{\bf Remark.} The proofs  of the theorem and the proposition, given in this paper, 
are still  valid  for more than three vector spaces when the format is the boundary
format. 

\smallskip

{\small I would like to thank  L.Gruson, M.Meulien and N.Perrin for their
help.}

\end{document}